\documentclass[12pt]{amsart}
\usepackage{amsmath,amssymb,amsbsy,amsfonts,amsthm,latexsym,amsopn,amstext,amsxtra,euscript,
amscd}
\usepackage{hyperref}
\usepackage{color}
\begin{document}

\title[Properties and applications of the PDF ...]{Properties and applications\\
of the prime detecting function:\\
infinitude of twin primes,\\
asymptotic law of distribution\\
of prime pairs\\
differing by an even number}

\author[R. M.~Abrarov]{R. M.~Abrarov$^1$}
\email{$^1$rabrarov@gmail.com}
\author[S. M.~Abrarov]{S. M.~Abrarov$^2$}
\email{$^2$absanj@gmail.com}

\date{September 29, 2011}

\begin{abstract}

The prime detecting function (PDF) approach can be effective instrument in the investigation of numbers. The PDF is constructed by recurrence sequence - each successive prime adds a sieving factor in the form of PDF. With built-in prime sieving features and properties such as simplicity, integro-differentiability and configurable capability for a wide variety of problems, the application of PDF leads to new interesting results. As an example, in this exposition we present proofs of the infinitude of twin primes and the first Hardy-Littlewood conjecture for prime pairs (the twin prime number theorem). On this example one can see that application of PDF is especially effective in investigation of asymptotic problems in combination with the proposed method of test and probe functions.
\end{abstract}
\maketitle
\section{\textbf{Properties of the prime detecting function}}
In our previous paper \cite{RAbrarov} we introduced the prime detecting function (PDF) with outcomes 1 for primes and 0 for composite numbers. In the general form the PDF looks like
\begin{equation} \label{Eq_1}
\begin{aligned}
p_{}^\delta \left( n \right): = & \left( {1 - \delta \left( {\frac{n}{{{p_1}}} - {p_1}} \right)} \right)\left( {1 - \delta \left( {\frac{n}{{{p_2}}} - {p_2}} \right)} \right)...\\
&\left( {1 - \delta \left( {\frac{n}{{{p_k}}} - {p_k}} \right)} \right)... \ , \qquad n \ge 2.
\end{aligned}
\end{equation}
where ${p_k} \ {\rm{ }}\left( {k = 1,{\rm{ 2, 3, }}...} \right)$ - all prime numbers in ascending order (throughout the article $k,{\rm{ }}l,{\rm{ }}m,{\rm{ }}n,{\rm{ }}q,{\rm{ }}r,{\rm{ }}s$
are integer numbers) and delta-function
\begin{equation} \label{Eq_2}
\delta\left(x\right):=
\begin{cases} 
1, \text{if $x \in \mathbb{N}_0=\{0, 1, 2,\dots\}$}  \\
0, \text{if $x \notin \mathbb{N}_0$}
\end{cases} .
\end{equation}
Directly from \eqref{Eq_1} and \eqref{Eq_2} follows if $n$ coincides with any ${p_k}{\rm{ }}\left( {k = 1,{\rm{ 2, 3, }}...} \right)$ then $p_{}^\delta \left( n \right) = 1$, else $p_{}^\delta \left( n \right) = 0$.

Consider in more detail the properties of this function. We do not know in advance the values of prime numbers ${p_k}$ for any $k$, but we can determine them successively by procedure described below. 

Initially, we assume that no single value of primes is known and the PDF has the form $p_0^\delta \left( n \right) = 1$. This function must be recharged timely by the known prime numbers ${p_k}$, and it will provide information about further numbers whether it is prime or composite. It is a recurrent routine - once the prime has been found, it induces change in the form of the PDF that, in turn, determines consecutively further primes.

Let us examine positive integers from the first integer following 1. It is 2. Consequently, we have $p_0^\delta \left( 2 \right) = 1$. It indicates that we found the fist prime number ${p_1} = 2$. Next, we have to charge it into PDF. After that it takes the form $p_1^\delta \left( n \right) = \left( {1 - \delta \left( {\frac{n}{2} - 2} \right)} \right)$. The next integer is 3 and $p_1^\delta \left( 3 \right) = \left( {1 - \delta \left( {\frac{3}{2} - 2} \right)} \right) = 1$.
Our function shows correct result that the second prime number ${p_2} = 3$. Again, it is necessary to charge prime ${p_2}$ in the function $p_1^\delta \left( n \right)$ getting the new form\\
$p_2^\delta \left( n \right) = \left( {1 - \delta \left( {\frac{n}{2} - 2} \right)} \right)\left( {1 - \delta \left( {\frac{n}{3} - 3} \right)} \right)$. 
Next natural number is 4 and $p_2^\delta \left( 4 \right) = \left( {1 - \delta \left( {\frac{4}{2} - 2} \right)} \right)\left( {1 - \delta \left( {\frac{4}{3} - 3} \right)} \right)$. The first factor in the product gives zero and, correspondingly, $p_2^\delta \left( 4 \right) = 0$. It means, 4 is composite number inducing no change in the form of the PDF. Further, we have $n = 5$. None of the factors in the product is zero. Therefore $p_2^\delta \left( 5 \right) = 1$ and ${p_3} = 5$. We have to charge this prime in our function $p_2^\delta \left( n \right)$ and now we get
$p_3^\delta \left( n \right) = \left( {1 - \delta \left( {\frac{n}{2} - 2} \right)} \right)\left( {1 - \delta \left( {\frac{n}{3} - 3} \right)} \right)\left( {1 - \delta \left( {\frac{n}{5} - 5} \right)} \right)$.

And so on, going consistently over all integer numbers, we can find out what the numbers are - prime or composite. Since there are infinitely many primes (Euclid's theorem \cite{Goldston, Hardy, Ribenboim, RAbrarov}), the PDF ultimately becomes
\begin{equation} \label{Eq_3}
\begin{aligned}
 {p^\delta }\left( n \right) = p_\infty ^\delta \left( n \right) = & \left( {1 - \delta \left( {\frac{n}{2} - 2} \right)} \right)\left( {1 - \delta \left( {\frac{n}{3} - 3} \right)} \right)\\
& \left( {1 - \delta \left( {\frac{n}{5} - 5} \right)} \right)... \\ 
= & \prod\limits_{p \ {\rm{primes}}} {\left( {1 - \delta \left( {\frac{n}{p} - p} \right)} \right)}, \qquad n \ge 2. \\
\end{aligned}
\end{equation}
The prime detecting function $p_k^\delta \left( n \right)$ can be defined by recurrence relations for $n = {\rm{2, 3, 4, }}...$ as
\begin{equation} \label{Eq_4}
\begin{cases}
   p_0^\delta \left( n \right) = 1,\\
   \text{if $p_{k - 1}^\delta \left( n \right) = 1 \Rightarrow {p_k} = n$ and replace $p_{k - 1}^\delta \left( n \right)$ with}\\
   \text{\qquad \qquad \qquad \quad \ \ \,$p_k^\delta \left( n \right) = p_{k - 1}^\delta \left( n \right)\left( {1 - \delta \left( {\frac{n}{{{p_k}}} - {p_k}} \right)} \right)$,}
\end{cases}
\end{equation}
Each ${p_k}$ is always a prime, since by recurrence procedure each number is checked on divisibility by each smaller prime and $p_{k - 1}^\delta \left( n \right) = 1$ means that ${p_k} = n$ does not have any smaller prime as a divisor. 

This function looks very simple and obvious. Nevertheless, it is significant enough to get important information about the numbers. Let us emphasize some important characteristics of the PDF. It is not difficult to observe that PDF utilizes an \emph{inclusion - exclusion principle}  \cite{Wiki-Inclusion}. Each bracket $\left( {1 - \delta \left( {\frac{n}{p} - p} \right)} \right)$ is a \emph{building block} of the PDF acting as \emph{numerical (digital) filter or sieve}, which filters out all the numbers divisible by $p$ starting with the \emph{square} of this number, i.e. \emph{argument of each delta-function has a quadratic zero-point boundary}. It is an effective \emph{functional alternative} to the well-known sieve of Eratosthenes \cite{Halberstam, Hardy, Ribenboim}. Depending on the assigned problem we can construct the PDF \emph{in various blocks (brackets)} - not necessarily with recurrence procedure \eqref{Eq_4} - and, instead of a sequence of positive integers $n$, we can examine some \emph{special sets of numbers or some function}: $\left( {1 - \delta \left( {\frac{{f(n)}}{p} - p} \right)} \right)$. It means the PDF possesses substantial \emph{flexibility for applications to a wide variety of problems}. As we can see further due to \emph{simple properties of the delta-functions}, the PDF is \emph{especially effective in the investigation of asymptotic problems}. Quite natural generalization of PDF to the function of real argument leads to \emph{integro-differentiability of the function} and in this case, as we can see below, it becomes in explicit form a \emph{derivative from the prime counting function} \cite{Goldston, Hardy, Ribenboim}.

If we charge PDF only with first $k$ primes ${p_k}$
\begin{equation} \label{Eq_5}
\begin{aligned}
p_k^\delta \left( n \right) = & \left( {1 - \delta \left( {\frac{n}{{{p_1}}} - {p_1}} \right)} \right)\left( {1 - \delta \left( {\frac{n}{{{p_2}}} - {p_2}} \right)} \right)...\\
& \left( {1 - \delta \left( {\frac{n}{{{p_k}}} - {p_k}} \right)} \right) , \quad    n \ge 2
\end{aligned}
\end{equation} 
it can give us answer about testing number $n$ whether it is prime or not until $p_{k + 1}^2 - 1$. After that, PDF's outcome 0 indicates that testing number is composite (and certainly one of ${p_k}$ is divisor) while outcome 1 shows that number is relatively prime to first $k$ primes (none of ${p_k}$ is divisor), but not necessarily the prime number.

Note that for the prime counting function (the quantity of primes $p \le n$) we have
\begin{equation} \label{Eq_6}
	\pi \left( n \right) = \sum\limits_{m = 2}^n {{p^\delta }\left( m \right)} .
\end{equation}
Consequently for each $p_k^\delta \left( n \right)$, we can define corresponding summatory function
\begin{equation} \label{Eq_7}
	{\pi _k}\left( n \right) = \sum\limits_{m = 2}^n {p_k^\delta \left( m \right)}  = \sum\limits_{m = 2}^n {\prod\limits_{l = 1}^k {\left( {1 - \delta \left( {\frac{m}{{{p_l}}} - {p_l}} \right)} \right)} } ,
\end{equation}
which is counting quantity of numbers up to $n$ relatively prime to first $k$ primes ${p_k}$. Also we have
\begin{equation} \label{Eq_8}
 \pi \left( n \right) = \sum\limits_{m = 2}^n {{p^\delta }\left( m \right)}  = \sum\limits_{m = 2}^n {p_\infty ^\delta \left( m \right)}  = {\pi _\infty }\left( n \right).
\end{equation}

Further we investigate asymptotic densities of these functions. Obviously, for 
$$\pi_1 \left(n\right)=\sum _{m=2}^{n}p_{{1}}^{\delta} \left(m\right)=\sum _{\text{odd numbers $\leq n$}} 1$$
we have an asymptotic density of odd numbers
\begin{equation} \label{Eq_9}
\mathop{\lim }\limits_{n\to \infty } \frac{\pi _1 \left(n\right)}{n}=\frac{1}{2}=\left(1-\frac{1}{p_{1} } \right) . 
\end{equation}
Using the following two remarkable properties of the delta-functions (further $\left\lfloor {\ } \right\rfloor {\rm{ and }}\left\lceil {\ } \right\rceil $ are floor and ceiling functions, respectively, and $\left\{ {\ } \right\}$ is fractional part of the number \cite{Wiki-Floor})
\begin{equation} \label{Eq_10}
\begin{aligned}
 & \delta \left(\frac{n }{p_{m_1} } -p_{m_1} \right)\dots  \delta \left(\frac{n}{p_{m_{i}}} - p_{m_{i}} \right) \delta \left(\frac{n}{p_{k}} - p_{k} \right) = \\
 & \delta \left( {\frac{n}{p_{m_1} }} \right) \dots \delta \left( {\frac{n}{p_{m_{i}} }} \right) \delta \left(\frac{n}{p_{k}} - p_{k} \right) = \\ 
 & \delta \left( {\frac{n}{{{p_{m_1}} \dots {p_{m_{i}}} \cdot p_{k} }} - \left\lceil {\frac{{{p_k}}}{{{p_{m_1}} \dots {p_{m_{i}}}}}} \right\rceil } \right), \quad m_1< \ldots <m_i<k,
\end{aligned} 
\end{equation}
and
\begin{equation} \label{Eq_11}
\begin{aligned}
& \mathop {\lim}\limits_{n \to \infty } \frac{{\sum\limits_{m = 2}^n {\delta \left( {\frac{m}{p_{k} \cdot q} - r_k} \right)} }}{n} = \frac{1}{p_{k} \cdot q} = \\
& \frac{1}{p_{k}} \mathop {\lim}\limits_{n \to \infty } \frac{{\sum\limits_{m = 2}^n {\delta \left( {\frac{m}{q} - r_q} \right)} }}{n} \text{     (for any finite $r_q$) ,}
\end{aligned}
\end{equation}
it is not difficult to get recurrence relation (assuming that $\mathop{\lim }\limits_{n\to \infty } \frac{\pi _{k-1} \left(n\right)}{n}$ exists)
\begin{equation} \label{Eq_12}
\begin{aligned}
& \mathop{\lim }\limits_{n\to \infty } \frac{\pi _{k} \left(n\right)}{n} =\\
 & \mathop{\lim }\limits_{n\to \infty } \frac { \sum\limits_{m = 2}^n {\left( {1 - \delta \left( {\frac{m}{{{p_k}}} - {p_k}} \right)} \right)\prod\limits_{l = 1}^{k-1} {\left( {1 - \delta \left( {\frac{m}{{{p_l}}} - {p_l}} \right)} \right)} } }{n} =\\
& \mathop{\lim }\limits_{n\to \infty } \frac{\pi _{k-1} \left(n\right)}{n}-\mathop{\lim }\limits_{n\to \infty } \frac { \sum\limits_{m = 2}^n {{ \delta \left( {\frac{m}{{{p_k}}} - {p_k}} \right)} \prod\limits_{l = 1}^{k-1} {\left( {1 - \delta \left( {\frac{m}{{{p_l}}} - {p_l}} \right)} \right)} } }{n} = \\
& \left( {1 - \frac{1}{{{p_k}}}} \right)\mathop{\lim }\limits_{n\to \infty } \frac{\pi _{k-1} \left(n\right)}{n} .
\end{aligned}
\end{equation}
From \eqref{Eq_9}, repeatedly applying the recurrence relation \eqref{Eq_12}, we obtain
\small
\begin{equation} \label{Eq_13}
\mathop {\lim}\limits_{n \to \infty } \frac{{{\pi _k}\left( n \right)}}{n} = \left( {1 - \frac{1}{{{p_1}}}} \right)\left( {1 - \frac{1}{{{p_2}}}} \right) \cdot ... \cdot \left( {1 - \frac{1}{{{p_k}}}} \right) = \prod\limits_{l = 1}^k {\left( {1 - \frac{1}{{{p_l}}}} \right)} {\rm{   }}{\rm{.}}
\end{equation}
\normalsize
which is valid for any $k \in \{1, 2, 3,\ldots\}$ .
In the limit of infinitely large $k$ in \eqref{Eq_13} from \eqref{Eq_5} and \eqref{Eq_8}, we conclude that asymptotic density of primes is
\begin{equation} \label{Eq_14}
\mathop{\lim }\limits_{n\to \infty } \frac{\pi _{\infty} \left(n\right)}{n} = \mathop{\lim }\limits_{n\to \infty } \frac{\pi \left(n\right)}{n} = \prod\limits_{\text{$p$ prime}} {\left( {1 - \frac{1}{p}} \right)} . 
\end{equation}

	We can generalize argument $n$ of PDF for any real number $x$ simply accepting
\begin{equation} \label{Eq_15}
	{\hat p^\delta }\left( x \right): = {p^\delta }\left( {\left\lceil x \right\rceil } \right).
\end{equation}
For this function there is corresponding analog of the prime counting function \eqref{Eq_6}
\begin{equation} \label{Eq_16}
\hat \pi \left( x \right): = \pi \left( {\left\lfloor x \right\rfloor } \right) + \left( {\pi \left( {\left\lfloor {x + 1} \right\rfloor } \right) - \pi \left( {\left\lfloor x \right\rfloor } \right)} \right)\left\{ x \right\},
\end{equation}
which is represented graphically as the values of the original prime counting function for integer arguments $\pi \left( {\left\lfloor x \right\rfloor } \right)$ \eqref{Eq_6} connected by the linear function $\left( {\pi \left( {\left\lfloor {x + 1} \right\rfloor } \right) - \pi \left( {\left\lfloor x \right\rfloor } \right)} \right)\left\{ x \right\}$. Then for the left derivative of \eqref{Eq_16} we get
\begin{equation} \label{Eq_17}
	\frac{{{d_{left}}\hat \pi \left( x \right)}}{{dx}} = {\hat p^\delta }\left( x \right)
\end{equation}
and
\begin{equation} \label{Eq_18}
	\hat \pi \left( x \right) = \int_2^x {{{\hat p}^\delta }\left( y \right)dy}.
\end{equation}
	If we want to see the prime counting function in a traditional way as a unit step function with jumps at the primes then for the corresponding PDF we have to perform few steps. First, we generalize \eqref{Eq_1} for any real $x$
\small
\begin{equation} \label{Eq_19}
\begin{aligned}
 p_{}^\delta \left( x \right): = & \left( {\delta \left( {x - 2} \right) - \delta \left( {\frac{x}{{{p_1}}} - {p_1}} \right)} \right)\left( {\delta \left( {x - 2} \right) - \delta \left( {\frac{x}{{{p_2}}} - {p_2}} \right)} \right)... \\ 
& {\rm{              }}\left( {\delta \left( {x - 2} \right) - \delta \left( {\frac{x}{{{p_k}}} - {p_k}} \right)} \right) \ldots \quad . 
\end{aligned}
\end{equation}
\normalsize
The function above is zero everywhere on real $x$ except primes on which the function is 1. Define the PDF for real argument
\begin{equation} \label{Eq_20}
	p_{}^{{\delta ^D}}\left( x \right): = p_{}^\delta \left( x \right){\delta ^D}\left( {x - \delta \left( x \right)x} \right) ,
\end{equation}
where ${\delta ^D}\left( x \right)$  -  Dirac delta-function (distribution) \cite{Estrada}, \cite{SAbrarov}\footnote{In our work \cite{SAbrarov} we showed simple examples involving number series, their reciprocals and some identities how the Dirac delta can be introduced from the new point of view. Compared to the Stieltjes integration methods widely utilized in number theory, the Dirac delta function approach turns out to be more efficient, natural and simple for understanding especially with the stepwise/discrete functions due to its integro-differentiability and other remarkable properties.}, \cite{Wiki-Dirac} with the sifting property
$$\int_{ - y - \varepsilon }^{ + y + \varepsilon } {f\left( x \right){\delta ^D}\left( {x - y} \right)dx}  = f\left( y \right){\rm{,  }} \quad \forall \varepsilon  > 0.$$
It should be noted that ${\delta ^D}\left( {x - \delta \left( x \right)x} \right) = \sum\limits_{m = 0}^\infty  {{\delta ^D}\left( {x - m} \right)}$ is nothing but a Dirac comb \cite{Wiki-Dirac} for $x \ge 0$, while the PDF $p_{}^{{\delta ^D}}\left( x \right)$ is linear combination of the delta distributions with spikes on primes representing a "reminder" left after sampling of the comb. It means, $p_{}^{{\delta ^D}}\left( x \right)$ has all the properties of a delta distribution (integro-differentiability, well defined Fourier, Laplace and other transforms, etc.). 
The prime counting function will be
\begin{equation} \label{Eq_21}
	\pi \left( x \right) = \int_{ - \infty }^x {{p^{{\delta ^D}}}\left( y \right)dy}  = \int_2^x {{p^{{\delta ^D}}}\left( y \right)dy},
\end{equation}
and correspondingly
\begin{equation} \label{Eq_22}
	\frac{{d\pi \left( x \right)}}{{dx}} = {p^{{\delta ^D}}}\left( x \right).
\end{equation}
We see that the PDF for real argument is in explicit form just derivative from the prime counting function.

\section{\textbf{Infinitude of twin primes,\\ the asymptotic law of distribution of prime pairs\\ differing by an even number}}

	Twin primes (or prime twins) are pairs of prime numbers that only differ by two. Some \ examples of twin primes \ (3, 5), (5, 7), (11, 13), (17, 19), (29, 31), \dots . There is only one pair of primes that differ less by two. It is pair (2, 3), which we can call the unique twin primes and put them apart from the set of twin primes. Obviously, all \ twin primes \ are \ pairs of \ successive \ odd \ numbers \ and \ except of \ (3, 5), (5, 7) there is no more two twin primes consisting of three successive odd numbers as in any such sequence bigger than 3, certainly one number is not a prime (must be divisible by 3). 
	There are many heuristic and computational \cite{Caldwell, Goldston, Hardy23, Ribenboim} evidences that the twin primes are infinitely many (twin primes conjecture), but infinitude of twin primes has not been proved until now. There have been several outstanding results towards the twin primes conjecture. In 1919 Brun showed that the sum of reciprocals of the twin primes (unlike the sum of reciprocals of primes) is either finite or convergent and comes to some constant $B = 1.9021604\dots,$ called Brun's constant
\begin{equation} \label{Eq_23}
 \left( {\frac{1}{3} + \frac{1}{5}} \right) + \left( {\frac{1}{5} + \frac{1}{7}} \right) + \left( {\frac{1}{{11}} + \frac{1}{{13}}} \right) + \left( {\frac{1}{{17}} + \frac{1}{{19}}} \right) + ... = B.
\end{equation}
Convergence signifies that the twin primes are quite sparse among the sequence of all primes. In 1920 Brun also gave upper bound for quantity of twin primes up to $x$, which is at present time improved to
\begin{equation} \label{Eq_24}
 \le 3.42 \cdot 2\prod\limits_{p > 2} {\left( {1 - \frac{1}{{{{\left( {p - 1} \right)}^2}}}} \right)} \frac{x}{{{{\left( {\log x} \right)}^2}}} .	\end{equation}
In famous paper \cite{Hardy23} (1923) Hardy and Littlewood on heuristic grounds conjectured that the estimation for quantity of twin primes up to $x$ should be the following:
\begin{equation} \label{Eq_25}
 	2{C^{twin}}\int_2^x {\frac{{dx}}{{{{\left( {\log x} \right)}^2}}}}
\end{equation}
with the twin prime constant (also for constant there are designations $C_2$ and $\Pi _2$)
\begin{equation} \label{Eq_26}
\begin{aligned}
{C^{twin}}  & = \prod\limits_{p > 2} {\frac{{p\left( {p - 2} \right)}}{{{{\left( {p - 1} \right)}^2}}} } = \prod\limits_{p > 2} {\left( {1 - \frac{1}{{{{\left( {p - 1} \right)}^2}}}} \right) }\\
& = 0.66016181584686957392781211\dots  .
\end{aligned}
\end{equation}
Chen proved that for infinitely many primes $p$, the number $p + 2$ is at most product of two primes \cite{Chen, Halberstam}. Quite recently another spectacular progress on this problem has been obtained \cite{Goldston-Pintz, Soundararajan}: team of mathematicians Goldston, Pintz and Yildirim proved that there are infinitely many consecutive primes, which are closer than any arbitrarily small multiple of the average spacing between primes
\begin{equation} \label{Eq_27}
 \Delta : = \mathop {\lim \inf }\limits_{n \to \infty } \frac{{{p_{n + 1}} - {p_n}}}{{\log n}} = 0.
\end{equation}
Assuming a very regular distribution of primes in arithmetic progressions same team showed that there are infinitely many pairs of primes differing by 16 or less. The current state of the subject can be found in Refs. \cite{Goldston, Goldston-Pintz, Korevaar, Soundararajan, Yildirim}.

The PDF approach can help to get further advance on the problem. By analogy with the PDF for any $n \ge 3$ we can construct a twin primes detecting function for pair $\left( {n,n + 2} \right)$ with outcomes 1 for twin primes, and 0 otherwise. To do this we apply the product of two PDFs
\begin{equation} \label{Eq_28}
	{p^{\delta :twin}}\left( n \right) = {p^\delta }\left( n \right){p^\delta }\left( {n + 2} \right),{\rm{      }} \quad n \ge 3.
\end{equation}
It is clear that such a function gives for $n$ a value 1 if  and only if ${p^\delta }\left( n \right)$ and ${p^\delta }\left( {n + 2} \right)$ both are 1, i.e. pair $\left( {n,n + 2} \right)$ are prime numbers (twin primes).
Recalling product formula \eqref{Eq_3} for each ${p^\delta }\left( n \right)$ and ${p^\delta }\left( {n + 2} \right)$ we can represent ($n \ge 3$)
\begin{equation} \label{Eq_29}
\begin{aligned}
 {p^{\delta :twin}}\left( n \right) = & \left( {1 - \delta \left( {\frac{n}{2} - 2} \right)} \right)\left( {1 - \delta \left( {\frac{{n + 2}}{2} - 2} \right)} \right)\\
& \left( {1 - \delta \left( {\frac{n}{3} - 3} \right)} \right)\left( {1 - \delta \left( {\frac{{n + 2}}{3} - 3} \right)} \right) \\
& \left( {1 - \delta \left( {\frac{n}{5} - 5} \right)} \right)\left( {1 - \delta \left( {\frac{{n + 2}}{5} - 5} \right)} \right) \dots . 
\end{aligned}
\end{equation}
Consider two pairs of products: $$\left( {1 - \delta \left( {\frac{n}{2} - 2} \right)} \right)\left( {1 - \delta \left( {\frac{{n + 2}}{2} - 2} \right)} \right)$$ and $$\left( {1 - \delta \left( {\frac{n}{3} - 3} \right)} \right)\left( {1 - \delta \left( {\frac{{n + 2}}{3} - 3} \right)} \right).$$ 
We have
\small
\begin{equation} \label{Eq_30}
\begin{aligned}
 & \left( {1 - \delta \left( {\frac{n}{2} - 2} \right)} \right) \left( {1 - \delta \left( {\frac{{n + 2}}{2} - 2} \right)} \right) = 1 - \delta \left( {\frac{n}{2} - 2} \right) - \\
 & \delta \left( {\frac{{n + 2}}{2} - 2} \right) + \delta \left( {\frac{n}{2} - 2} \right)\delta \left( {\frac{{n + 2}}{2} - 2} \right). 
\end{aligned}
\end{equation}
\normalsize
Sum of two last terms $$ - \delta \left( {\frac{{n + 2}}{2} - 2} \right) + \delta \left( {\frac{n}{2} - 2} \right)\delta \left( {\frac{{n + 2}}{2} - 2} \right)$$ for $n > 2$ are always zero because either $n$ is odd number and both summands are zero or $n$ is even and summands have opposite values $ - 1 + 1 = 0$. Thus, product is simplified
\begin{equation} \label{Eq_31}
	\left( {1 - \delta \left( {\frac{n}{2} - 2} \right)} \right)\left( {1 - \delta \left( {\frac{{n + 2}}{2} - 2} \right)} \right) = 1 - \delta \left( {\frac{n}{2} - 2} \right).
\end{equation}
In the second product
\small
\begin{equation} \label{Eq_32}
\begin{aligned}
& \left( {1 - \delta \left( {\frac{n}{3} - 3} \right)} \right) \left( {1 - \delta \left( {\frac{{n + 2}}{3} - 3} \right)} \right) = 1 - \delta \left( {\frac{n}{3} - 3} \right) - \\
& \delta \left( {\frac{{n + 2}}{3} - 3} \right) + \delta \left( {\frac{n}{3} - 3} \right)\delta \left( {\frac{{n + 2}}{3} - 3} \right) 
\end{aligned}
\end{equation}
\normalsize
term $\delta \left( {\frac{n}{3} - 3} \right)\delta \left( {\frac{{n + 2}}{3} - 3} \right)$ always zero because 3 never can divide both $n$ and $n + 2$. So we have
\footnotesize
\begin{equation} \label{Eq_33}
	\left( {1 - \delta \left( {\frac{n}{3} - 3} \right)} \right)\left( {1 - \delta \left( {\frac{{n + 2}}{3} - 3} \right)} \right) = 1 - \delta \left( {\frac{n}{3} - 3} \right) - \delta \left( {\frac{{n + 2}}{3} - 3} \right).
\end{equation}
\normalsize
Above simplification for prime 3 can be applied for any bigger prime numbers. Hence, the twin primes detecting function becomes
\footnotesize
\begin{equation} \label{Eq_34}
\begin{aligned}
 {p^{\delta :twin}}\left( n \right) = & \left( {1 - \delta \left( {\frac{n}{2} - 2} \right)} \right) \cdot \left( {1 - \delta \left( {\frac{n}{3} - 3} \right) - \delta \left( {\frac{{n + 2}}{3} - 3} \right)} \right) \cdot  \\ 
& \left( {1 - \delta \left( {\frac{n}{5} - 5} \right) - \delta \left( {\frac{{n + 2}}{5} - 5} \right)} \right) \cdot ...{\rm{  ,  }} \quad n \ge 3.
\end{aligned} 
\end{equation}
\normalsize
In general, the expression for the twin primes detecting function we represent as
\footnotesize
\begin{equation} \label{Eq_35}
\begin{aligned}
 {p^{\delta :twin}}\left( n \right) = & \left( {1 - \delta \left( {\frac{n}{{{p_1}}} - {p_1}} \right)} \right)\left( {1 - \delta \left( {\frac{n}{{{p_2}}} - {p_2}} \right) - \delta \left( {\frac{{n + 2}}{{{p_2}}} - {p_2}} \right)} \right) \dots \\ 
& \left( {1 - \delta \left( {\frac{n}{{{p_k}}} - {p_k}} \right) - \delta \left( {\frac{{n + 2}}{{{p_k}}} - {p_k}} \right)} \right) \dots , \quad n \ge 3, 
\end{aligned}
\end{equation}
\normalsize
where ${p_1} = 2$, ${p_2} = 3$, ${p_3} = 5$ and so on.

Now following the method discussed in previous section, we try to find asymptotic densities for twin pairs relatively prime to first $k$ prime numbers. The corresponding twin pairs detecting function is
\footnotesize
\begin{equation} \label{Eq_36}
\begin{aligned}
 p_k^{\delta :twin}\left( n \right) = & \left( {1 - \delta \left( {\frac{n}{{{p_1}}} - {p_1}} \right)} \right)\\
& \prod\limits_{m = 2}^k {\left( {1 - \delta \left( {\frac{n}{{{p_m}}} - {p_m}} \right) - \delta \left( {\frac{{n + 2}}{{{p_m}}} - {p_m}} \right)} \right)}, \quad n \ge 3,
\end{aligned}
\end{equation}
\normalsize
where ${p_m}{\rm{  }}\left( {m = 1,{\rm{ 2, 3, }}...{\rm{, }}k} \right)$ - first $k$ prime numbers in ascending order.
Corresponding summatory function
\begin{equation} \label{Eq_37}
	\pi _k^{twin}\left( n \right) = \sum\limits_{m = 3}^n {p_k^{\delta {\rm{:}}twin}\left( m \right)} \quad .
\end{equation}
For $\pi_1^{twin} \left(n\right)$ counting twin pairs relatively prime to 2, i.e. all odd pairs beginning from 3 up to $n$
$$\pi_1^{twin} \left(n\right)=\sum _{m=2}^{n}p_{{1}}^{\delta} \left(m\right)=\sum _{3\leq\text{odd numbers $\leq n$}} 1,$$
we have same asymptotic density as in \eqref{Eq_9}
\begin{equation} \label{Eq_38}
\mathop{\lim }\limits_{n\to \infty } \frac{\pi _1^{twin} \left(n\right)}{n}=\frac{1}{2}=\left(1-\frac{1}{p_{1} } \right) . 
\end{equation}
Unlike the PDF \eqref{Eq_5}, disclosure of the brackets in \eqref{Eq_36} leads also to the terms such as
\footnotesize
\begin{equation} \label{Eq_39}
\begin{aligned}
 & \delta \left(\frac{n }{p_{k_1} } -p_{k_1} \right)\dots \delta \left(\frac{n}{p_{k_{i}}} - p_{k_{i}} \right)  
   \delta \left(\frac{n+2}{p_{l_1} } -p_{l_1} \right)\dots \delta \left(\frac{n+2}{p_{l_{j}}} - p_{l_{j}} \right) = \\
 & \delta \left( \frac{n}{{q_1} }- \left\lceil \frac{p_{k_{i}}^2}{q_1} \right\rceil  \right)\delta \left( \frac{n+2}{{q_2} }- \left\lceil \frac{p_{l_{j}}^2}{q_2} \right\rceil  \right), \\
 & k_1< \ldots <k_{i-1}<k_i, \quad l_1< \ldots <l_{j-1}<l_j, \\
 & q_1=p_{k_1}\ldots p_{k_i}, \quad q_2=p_{l_1}\ldots p_{l_j}, \quad \left(q_1,\ q_2 \right)=1,
\end{aligned} 
\end{equation}
\normalsize
Above product of delta-functions can be reduced to a single delta-function of the form $\delta \left( {\frac{{m + s}}{q} - r} \right)$ with $q = q_1 \cdot q_2 $ since $q_1$ and $q_2$ are relatively prime, and we always can find finite integers $m_1$ and $m_2$ such as $q_1 \cdot m_1= q_2 \cdot m_2 +2$. So we rewrite \eqref{Eq_39} as
\small
\begin{equation} \label{Eq_40}
\begin{aligned}
 & \delta \left( \frac{n+q_1 m_1 }{{q_1}}- \left\lceil \frac{p_{k_{i}}^2}{q_1} + m_1 \right\rceil  \right) \delta \left( \frac{n+2 + q_2 m_2}{{q_2} }- \left\lceil \frac{p_{l_{j}}^2}{q_2} + m_2 \right\rceil  \right) = \\
 & \delta \left( \frac{n+s}{q} - r  \right), \quad s=q_1 m_1=q_2 m_2+2, \quad r \text{ - some finite integer}.
\end{aligned} 
\end{equation}
\normalsize
Applying \eqref{Eq_39}, \eqref{Eq_40} and analog of \eqref{Eq_11}
\begin{equation} \label{Eq_41}
\begin{aligned}
& \mathop {\lim}\limits_{n \to \infty } \frac{{\sum\limits_{m = 2}^n {\delta \left( {\frac{m+s_k}{p_{k} \cdot q} - r_k} \right)} }}{n} = \frac{1}{p_{k} \cdot q} = \\
& \frac{1}{p_{k}} \mathop {\lim}\limits_{n \to \infty } \frac{{\sum\limits_{m = 2}^n {\delta \left( {\frac{m+s_q}{q} - r_q} \right)} }}{n} \text{     (for any finite $s_q$ and $r_q$) ,}
\end{aligned}
\end{equation}
it is not difficult to get recurrence relation for $k\geq2$ (assuming that $\mathop{\lim }\limits_{n\to \infty } \frac{\pi _{k-1}^{twin} \left(n\right)}{n}$ exists)
\begin{equation} \label{Eq_42}
\begin{aligned}
& \mathop{\lim }\limits_{n\to \infty } \frac{\pi _{k}^{twin} \left(n\right)}{n} =\\
 & \mathop{\lim }\limits_{n\to \infty } \frac { \sum\limits_{m = 2}^n {\left( {1 - \delta \left( {\frac{m}{{{p_k}}} - {p_k}} \right) - \delta \left( {\frac{m+2}{{{p_k}}} - {p_k}} \right)} \right) p_{k-1}^{\delta :twin}\left( m \right) } }{n} =\\
& \mathop{\lim }\limits_{n\to \infty } \frac{\pi _{k-1}^{twin} \left(n\right)}{n}-\mathop{\lim }\limits_{n\to \infty } \frac { \sum\limits_{m = 2}^n {{ \delta \left( {\frac{m}{{{p_k}}} - {p_k}} \right)} p_{k-1}^{\delta :twin}\left( m \right) } }{n} - \\
& \mathop{\lim }\limits_{n\to \infty } \frac { \sum\limits_{m = 2}^n {{ \delta \left( {\frac{m+2}{{{p_k}}} - {p_k}} \right)} p_{k-1}^{\delta :twin}\left( m \right) } }{n} = \left( {1 - \frac{2}{{{p_k}}}} \right)\mathop{\lim }\limits_{n\to \infty } \frac{\pi _{k-1}^{twin} \left(n\right)}{n} .
\end{aligned}
\end{equation}

From \eqref{Eq_38}, applying recurrence relation \eqref{Eq_42} as many time as we want, we can obtain an asymptotic density
\begin{equation} \label{Eq_43}
\begin{aligned}
	\mathop {\lim}\limits_{n \to \infty } \frac{{\pi _k^{twin}\left( n \right)}}{n} & = \left( {1 - \frac{1}{{{p_1}}}} \right)\left( {1 - \frac{2}{{{p_2}}}} \right) \cdot ... \cdot \left( {1 - \frac{2}{{{p_k}}}} \right) \\ 
	& = \left( {1 - \frac{1}{{{p_1}}}} \right)\prod\limits_{m = 2}^k {\left( {1 - \frac{2}{{{p_m}}}} \right)},
\end{aligned}
\end{equation}
which is valid for any $k \in \{1,\ 2,\ 3,\ \ldots\}$ .
Further, we can infer that the twin prime counting function (quantity of twin primes $\left( {p,p + 2} \right)$ such that $p \le n$)
\begin{equation} \label{Eq_44}
	{\pi ^{twin}}\left( n \right) = \sum\limits_{m = 3}^n {{p^{\delta {\rm{:}}twin}}\left( n \right)}  = \sum\limits_{m = 3}^n {p_\infty ^{\delta {\rm{:}}twin}\left( n \right)}  = \pi _\infty ^{twin}\left( n \right)
\end{equation}
has asymptotic density
\small
\begin{equation} \label{Eq_45}
\begin{aligned}
	\mathop {\lim}\limits_{n \to \infty } \frac{{{\pi_\infty^{twin}}\left( n \right)}}{n} = \mathop {\lim}\limits_{n \to \infty } \frac{{{\pi ^{twin}}\left( n \right)}}{n} 
	& = \left( {1 - \frac{1}{2}} \right)\left( {1 - \frac{2}{3}} \right)\left( {1 - \frac{2}{5}} \right)\dots \\
	& = \left( {1 - \frac{1}{2}} \right)\prod\limits_{p > 2} {\left( {1 - \frac{2}{p}} \right)}.
\end{aligned}
\end{equation}
\normalsize

Now we combine the ratio of asymptotic density \eqref{Eq_43} and square of asymptotic density \eqref{Eq_13}
\begin{equation} \label{Eq_46}
	\mathop {\lim}\limits_{n \to \infty } \left( {\frac{{\pi _k^{twin}\left( n \right)}}{n}{{\left( {\frac{n}{{{\pi _k}\left( n \right)}}} \right)}^2}} \right)
\end{equation}
to get for each $k = 1,{\rm{ }}2,{\rm{ }}3,...$\\
\small
\leftline{$\qquad\quad
\begin{array}{l}
\begin{aligned}
\mathop {\lim}\limits_{n \to \infty } \left( {\frac{{\pi _1^{twin}\left( n \right)}}{n}{{\left( {\frac{n}{{{\pi _1}\left( n \right)}}} \right)}^2}} \right) & = \frac{1}{2}{\left( 2 \right)^2} = 2,\\
\mathop {\lim}\limits_{n \to \infty } \left( {\frac{{\pi _2^{twin}\left( n \right)}}{n}{{\left( {\frac{n}{{{\pi _2}\left( n \right)}}} \right)}^2}} \right) & = 2\left( {1 - \frac{2}{3}} \right){\left( {1 - \frac{1}{3}} \right)^{ - 2}} \\
 & = 2\left( {1 - \frac{1}{{{{\left( {3 - 1} \right)}^2}}}} \right), \\ 
 \mathop {\lim}\limits_{n \to \infty } \left( {\frac{{\pi _3^{twin}\left( n \right)}}{n}{{\left( {\frac{n}{{{\pi _3}\left( n \right)}}} \right)}^2}} \right) & = 2\left( {1 - \frac{1}{{{{\left( {3 - 1} \right)}^2}}}} \right)\left( {1 - \frac{2}{5}} \right){\left( {1 - \frac{1}{5}} \right)^{ - 2}}\\
 & = 2\left( {1 - \frac{1}{{{{\left( {3 - 1} \right)}^2}}}} \right)\left( {1 - \frac{1}{{{{\left( {5 - 1} \right)}^2}}}} \right),
\end{aligned} 
\end{array}$}\\
\leftline{\qquad\qquad\dots}
\begin{equation} \label{Eq_47}
\begin{aligned}
\mathop {\lim}\limits_{n \to \infty } & \left( {\frac{{\pi _k^{twin}\left( n \right)}}{n}{{\left( {\frac{n}{{{\pi _k}\left( n \right)}}} \right)}^2}} \right) =\\ 
2 & \left( {1 - \frac{1}{{{{\left( {3 - 1} \right)}^2}}}} \right)
 \left( {1 - \frac{1}{{{{\left( {5 - 1} \right)}^2}}}} \right)...\left( {1 - \frac{1}{{{{\left( {{p_k} - 1} \right)}^2}}}} \right).
\end{aligned}
\end{equation}
\normalsize
Ultimately, when $k$ tends to infinity right hand side of expression \eqref{Eq_47} has limit (see \eqref{Eq_26})
\begin{equation} \label{Eq_48}
	2\prod\limits_{k = 2}^\infty  {\left( {1 - \frac{1}{{{{\left( {{p_k} - 1} \right)}^2}}}} \right)}  = 2{C^{twin}}.
\end{equation}

With tendency of $k$ to bigger numbers we see the process of "refining" of primes and twin primes among the integers relatively prime to first $k$ primes. Also the least composite number relatively prime to first $k$ primes is drifting to infinity, when $k \to \infty$.  For extremely large $k$ expanding to infinity the asymptotic of $\frac{{\pi _k^{twin}\left( n \right)}}{n}{\left( {\frac{n}{{{\pi _k}\left( n \right)}}} \right)^2}$ remains almost same (more and more close to constant $2{C^{twin}}$), i.e. for extremely large $k$ it is weakly sensitive for switching from one $k$ to another, signifying that the greatest contribution in asymptotic is getting from "pure" primes and "pure" twin primes only. These two processes of "refining" and approaching to constant $2{C^{twin}}$ are going parallel.  Ultimately, in the relation of asymptotic densities $\frac{{\pi _k^{twin}\left( n \right)}}{n}$ and ${\left( {\frac{{{\pi _k}\left( n \right)}}{n}} \right)^2}$ we get the ratio of only prime asymptotic densities  $\frac{{\pi _{}^{twin}\left( n \right)}}{n}$ and ${\left( {\frac{{\pi \left( n \right)}}{n}} \right)^2}$ (see \eqref{Eq_45} and \eqref{Eq_14})
\begin{equation} \label{Eq_49}
\begin{aligned}
	& \mathop {\lim}\limits_{n \to \infty } \left( {\frac{{\pi _\infty ^{twin}\left( n \right)}}{n}{{\left( {\frac{n}{{{\pi _\infty }\left( n \right)}}} \right)}^2}} \right) = \\
	& \mathop {\lim}\limits_{n \to \infty } \left( {\frac{{\pi _{}^{twin}\left( n \right)}}{n}{{\left( {\frac{n}{{\pi \left( n \right)}}} \right)}^2}} \right) = 2{C^{twin}}.
\end{aligned}
\end{equation}
Using the Prime Number Theorem \cite{Hardy, RAbrarov} stating that
\begin{equation} \label{Eq_50}
	\mathop {\lim}\limits_{n \to \infty } \left( {\frac{{\pi \left( n \right)}}{n}\log n} \right) = 1  \text{\ or \ }     \frac{{\pi \left( n \right)}}{n} \thicksim \frac{1}{{\log n}}
\end{equation}
 we can rewrite \eqref{Eq_49}
\begin{equation} \label{Eq_51}
\begin{aligned}
 	\pi _{}^{twin}\left( n \right) & \thicksim 2{C^{twin}}n{\left( {\frac{{\pi \left( n \right)}}{n}} \right)^2} \\ 
 	& \thicksim 2{C^{twin}}\frac{n}{{{{\left( {\log n} \right)}^2}}} \thicksim 2{C^{twin}}\int_2^n {\frac{{dx}}{{{{\left( {\log x} \right)}^2}}}} .
\end{aligned}
\end{equation}
Expression \eqref{Eq_51} explicitly asserts the infinitude of twin primes and gives the asymptotic law of distribution for twin primes. Thus, twin primes conjecture and Hardy - Littlewood conjecture for twin primes \eqref{Eq_25} is proved.

In our method, we define an asymptotic of the $\frac{{\pi _{}^{twin}\left( n \right)}}{n}$ by comparing it with the asymptotic of the function ${\left( {\frac{{\pi \left( n \right)}}{n}} \right)^2}$. In other words, here $\frac{{\pi _{}^{twin}\left( n \right)}}{n}$ is test-function with uncertain behaviour in the asymptotic and ${\left( {\frac{{\pi \left( n \right)}}{n}} \right)^2}$ is probe-function with known asymptotic. The same idea was used in our previous paper \cite{RAbrarov}, where we presented our version of proof of the Prime Number Theorem (in that work the test-function was $\frac{{\pi \left( n \right)}}{n}$ and probe-function was Harmonic Number $H\left( n \right)$).

Result obtained for twin primes can be readily generalized for any even distance $2k$ between prime pairs $\left( {p,p + 2k} \right)$ and for corresponding $2k$ - pair prime detecting function
\begin{equation} \label{Eq_52}
{p^{\delta :2k}}\left( n \right) = {p^\delta }\left( n \right){p^\delta }\left( {n + 2k} \right), \quad n \ge 3 .
\end{equation}
For this purpose we recall all inferences \eqref{Eq_28}-\eqref{Eq_45} and just examine that in
\[
\begin{aligned}
 & \left( {1 - \delta \left( {\frac{n}{p} - p} \right)} \right)\left( {1 - \delta \left( {\frac{{n + 2k}}{p} - p} \right)} \right) = 1
 - \delta \left( {\frac{n}{p} - p} \right) \\ 
 & - \delta \left( {\frac{{n + 2k}}{p} - p} \right) + \delta \left( {\frac{n}{p} - p} \right)\delta \left( {\frac{{n + 2k}}{p} - p} \right)
\end{aligned}
\]
four terms for any prime $p > 2$ dividing $k$ can be reduced to only two terms for $n \ge {p^2}$
\begin{equation} \label{Eq_53}
	\left( {1 - \delta \left( {\frac{n}{p} - p} \right)} \right)\left( {1 - \delta \left( {\frac{{n + k}}{p} - p} \right)} \right) = 1 - \delta \left( {\frac{n}{p} - p} \right)
\end{equation}
due to same reasons as we reduced in \eqref{Eq_30} and \eqref{Eq_31} for prime 2. It means that each prime $p$ dividing $2k$ contributes to asymptotic densities with factor $\left( {1 - \frac{1}{p}} \right)$ instead of factor $\left( {1 - \frac{2}{p}} \right)$. Consequently, in asymptotic densities for prime pairs differing by $2k$ we have to multiply on factors
\begin{equation} \label{Eq_54}
	\prod\limits_{\scriptstyle p > 2 \hfill \atop 
  \scriptstyle p\left| k \right. \hfill} {\frac{{1 - \frac{1}{p}}}{{1 - \frac{2}{p}}}}  = \prod\limits_{\scriptstyle p > 2 \hfill \atop 
  \scriptstyle p\left| k \right. \hfill} {\frac{{p - 1}}{{p - 2}}} .
\end{equation}
Finally, for the prime pairs counting function (quantity of prime pairs $\left( {p,p + 2k} \right)$ $p \le n$)  we have
\begin{equation} \label{Eq_55}
\begin{aligned}
	\pi _{}^{2k}\left( n \right) & \thicksim 2{C^{twin}}\prod\limits_{\scriptstyle p > 2 \hfill \atop 
  \scriptstyle p\left| k \right. \hfill} {\frac{{p - 1}}{{p - 2}}} \frac{n}{{{{\left( {\log n} \right)}^2}}} \\
  & \thicksim 2{C^{twin}}\prod\limits_{\scriptstyle p > 2 \hfill \atop 
  \scriptstyle p\left| k \right. \hfill} {\frac{{p - 1}}{{p - 2}}} \int_2^n {\frac{{dx}}{{{{\left( {\log x} \right)}^2}}}}  .
\end{aligned}
\end{equation}
It is what Hardy and Littlewood conjectured in 1923 about prime pairs (conjecture B) [\cite{Hardy23}, p.42]. Hence, the Hardy - Littlewood conjecture B now is proved. Obviously, formula \eqref{Eq_55} also claims a weaker statement - for every natural number $k$ there exist infinitely many pair of primes $p,{\rm{ }}q$ such that $q - p = 2k$ (without requiring them to be consecutive).

Remarkably, the prime detecting function approach and method with asymptotic densities for test and probe functions also can be applied to a wide variety of problems with sequences of prime polynomials (Goldbach conjecture, primes in arithmetic progressions, clusters of twin primes, $k - tuple$ conjecture, Cunningham chains, \dots) and  non-polynomial primes (such as Mersenne and Fermat primes) and other problems \cite{Caldwell, Ribenboim}.


\begin{thebibliography}{17}
\smallskip
\bibitem{RAbrarov}
R. M. Abrarov and S. M. Abrarov, \textit{Sieve procedure for the M\"obius prime-functions, the infinitude of primes and the prime number theorem}, 
\indent\verb+http://arxiv.org/abs/1004.1563+
\smallskip
\bibitem{SAbrarov}
S. M. Abrarov and R. M. Abrarov, \textit{Identities for number series and their reciprocals: Dirac delta-function approach},\\ \indent\verb+http://arxiv.org/pdf/0704.1936v3+
\smallskip
\bibitem{Caldwell}
C. K. Caldwell, \leftline{\textit{An amazing prime heuristic},} \indent\verb+http://www.utm.edu/staff/caldwell/preprints/Heuristics.pdf+
\smallskip
\bibitem{Chen}
J. R. Chen, \textit{On the representation of a larger even integer as the sum of a prime and the product of at most two primes}, Sci. Sinica \textbf{16} (1973) 157-176.
\smallskip
\bibitem{Estrada}
R. Estrada, R. P. Kanwal, \textit{A distributional approach to asymptotics. Theory and applications}, Second edition, Birkh\"auser, Boston, 2002.
\smallskip
\bibitem{Goldston}
D. A. Goldston, \leftline{\textit{Are there infinitely many twin primes?},} \indent\verb+http://arxiv.org/pdf/0710.2123v1+
\smallskip
\bibitem{Goldston-Pintz}
D. A. Goldston, J. Pintz and C. Y. Yildirim, \textit{Primes in tuples. I.}, Ann. of Math.\\ \textbf{170}, 2 (2009), 819-862, \indent\verb+http://dx.doi.org/10.4007/annals.2009.170.819+
Earlier version \indent\verb+http://arxiv.org/pdf/math/0508185v1+
\smallskip
\bibitem{Halberstam}
H. Halberstam and H.-E. Richert, \textit{Sieve Methods}, Academic Press, New York, 1974.
\smallskip
\bibitem{Hardy23}
G. H. Hardy and J. E. Littlewood, \textit{Some problems of 'Partitio Numerorum'; III: On the expression of a number as a sum of primes}, Acta Math. \textbf{44} (1923), 1 - 70; Reprinted in \textit{"Collected Papers of G. H. Hardy"}, Vol. I, pp. 561-630, Clarendon Press, Oxford, 1966.
\smallskip
\bibitem{Hardy}
G. H. Hardy and E. M. Wright, \textit{An introduction to the theory of numbers}, 5th ed., Oxford University Press, Oxford, 1979.
\smallskip
\bibitem{Korevaar}
J. Korevaar, \textit{Prime pairs and the zeta function}, Journal of Approximation Theory \ \textbf{158} \ (2009)\ 69-96, \indent\verb+http://dx.doi.org/10.1016/j.jat.2008.01.008+
\smallskip
\bibitem{Ribenboim}
P. Ribenboim, \textit{The New Book of Prime Number Records}, third edition, Springer-Verlag New York Inc., 1996.
\smallskip
\bibitem{Soundararajan}
K. Soundararajan, \textit{Small gaps between prime numbers: The work of Goldston-Pintz-Yildirim}, Bull. Amer. Math. Soc. (N.S.) \textbf{44} (2007) 1-18.
\smallskip
\bibitem{Wiki-Dirac}
Wikipedia contributors,
\href{http://en.wikipedia.org/wiki/Dirac_delta_function}{\textit{$"$Dirac delta functions$"$,}} Wikipedia, The Free Encyclopedia.
\smallskip
\bibitem{Wiki-Floor}
Wikipedia contributors,
\href{http://en.wikipedia.org/wiki/Floor_function}{\textit{$"$Floor and ceiling functions$"$,}} Wikipedia, The Free Encyclopedia.
\smallskip
\bibitem{Wiki-Inclusion}
Wikipedia contributors,
\href{http://en.wikipedia.org/wiki/Inclusion-exclusion_principle}{\textit{$"$Inclusion-exclusion principle$"$,}} Wikipedia, The Free Encyclopedia.
\smallskip
\bibitem{Yildirim}
C. Y. Yildirim, \textit{The Distribution of Primes: Conjectures vs. Hitherto Provables}, Further progress in analysis, World Sci. Publ., Hackensack, NJ, (2009) 75-108, \indent\verb+http://dx.doi.org/10.1142/9789812837332_0004+
\\
\end{thebibliography}
\end{document}